\documentclass{amsart}
\usepackage{amsmath}
\usepackage{amssymb} %\mathbb
\usepackage{amscd} %{\CD}:commutative diagram
\usepackage{graphicx} %attaching figure
\usepackage{epsfig}
\usepackage[dvips]{color}
\usepackage{slashbox}

\newtheorem{THM}{Theorem}[section]

\newtheorem{LEM}[THM]{Lemma}

\newtheorem{QUES}[THM]{Question}
\newtheorem{REM}[THM]{Remark}

\begin{document}

\title[Seifert fibered surgeries with distinct primitive/Seifert positions]
{Seifert fibered surgeries with distinct \\
primitive/Seifert positions} 

\author{Mario Eudave-Mu\~noz, \
Katura Miyazaki
and Kimihiko Motegi \\
}

\address{Instituto de Matematicas, Universidad Nacional Aut\'onoma de M\'exico, Circuito Exterior, Ciudad Universitaria 04510 M\'exico DF, Mexico}
\email{mario@matem.unam.mx}

\address{Faculty of Engineering, Tokyo Denki University, Tokyo 101--8457, 
Japan}
\email{miyazaki@cck.dendai.ac.jp}

\address{Department of Mathematics, Nihon University, 
Tokyo 156--8550, Japan}
\email{motegi@math.chs.nihon-u.ac.jp}

\date{}

\begin{abstract}
We call a pair $(K, m)$ of a knot $K$ in the $3$--sphere $S^3$ and an integer $m$ 
a Seifert fibered surgery if $m$--surgery on $K$ yields a Seifert fiber space. 
For most known Seifert fibered surgeries $(K, m)$, 
$K$ can be embedded in a genus $2$ Heegaard surface of $S^3$ 
in a primitive/Seifert position, the concept introduced by Dean
as a natural extension of primitive/primitive position defined by Berge. 
Recently Guntel has given an infinite family of
Seifert fibered surgeries each of which has distinct  
primitive/Seifert positions. 
In this paper we give yet other infinite families of Seifert fibered surgeries with 
distinct primitive/Seifert positions from a different point of view.
\end{abstract}

\maketitle
{%
\renewcommand{\thefootnote}{}
\footnotetext{2000 \textit{Mathematics Subject Classification.}
Primary 57M25, 57M50}
\footnotetext{ \textit{Key words and phrases.}
Dehn surgery, Seifert fiber space, primitive/Seifert position}
}%
 
\section{Introduction}
\label{sect:intro}

Let $(K, m)$ be a pair of a knot $K$ in $S^3$ and an integer $m$,
and denote by $K(m)$ the manifold obtained from $S^3$ by
$m$--surgery on $K$.
We say that $(K, m)$ is a \textit{Seifert fibered surgery}
if $K(m)$ is a Seifert fiber space.  
We regard that two Seifert fibered surgeries $(K, m)$ and $(K', m')$
are the \textit{same} if $K$ has the same knot type as $K'$
(i.e.\ $K$ is isotopic to $K'$ in $S^3$) and $m = m'$. 
For a genus $2$ handlebody $H$ and a simple closed curve $c$ in
$\partial H$, we denote $H$ with a $2$--handle attached along $c$
by $H[c]$.

Let $S^3 = V \cup_F W$ be a genus $2$ Heegaard splitting of $S^3$, 
i.e.\ $V$ and $W$ are genus $2$ handlebodies in $S^3$
with $V \cap W$ a genus $2$ Heegaard surface $F$. 
It is known that such a splitting is unique up to isotopy in $S^3$ \cite{Wa}. 
We say that a Seifert fibered surgery $(K, m)$ has a \textit{primitive/Seifert position} $(F, K', m)$ if 
$K'$ is a simple closed curve in a genus $2$ Heegaard surface $F$
such that $K' (\subset  S^3)$ has the same knot type as $K$ and
satisfies the following three conditions.

\begin{itemize}
\item 
$K'$ is \textit{primitive} with respect to $V$, 
i.e.\ $V[K']$ is a solid torus. 
\item
$K'$ is \textit{Seifert} with respect to $W$, 
i.e.\ $W[K']$ is a Seifert fiber space 
with the base orbifold $D^2(p, q)$ $(p, q \ge 2)$.  
\item
The \textit{surface slope} of $K'$ with respect to $F$ 
(i.e.\ the isotopy class in  $\partial N(K')$ represented by a component of $\partial N(K') \cap F$) 
coincides with the surgery slope $m$. 
\end{itemize}

For the primitive/Seifert position $(F, K', m)$ above,
we define the \textit{index set} $i(F, K', m)$
to be the set $\{ p, q \}$.

In general, 
if a knot $K$ in $S^3$ has a primitive/Seifert position with surface slope $m$, 
then $K$ is strongly invertible (\cite[Claim 5.3]{MM7}) and 
$K(m) \cong V[K] \cup W[K]$ is a Seifert fiber space or a connected sum of lens spaces. 
In particular, if $K$ is hyperbolic, then the latter case cannot happen by 
the positive solution to the cabling conjecture for strongly 
invertible knots \cite{EM}.  

The notion of primitive/Seifert position was introduced by Dean \cite{Dean} as a natural 
modification of Berge's primitive/primitive position \cite{Berge}. 
It is conjectured that all the lens surgeries have primitive/primitive positions \cite{Go90}. 
On the other hand, 
there are infinitely many Seifert fibered surgeries 
with no primitive/Seifert positions \cite{MMM, DMM, Tera}; 
nevertheless the majority of Seifert fibered surgeries have such positions.  
Let $(K, m)$ be a Seifert fibered surgery with 
two primitive/Seifert positions $(F_1, K_1, m)$ and $(F_2, K_2, m)$. 
Then, we say that $(F_1, K_1, m)$ and $(F_2, K_2, m)$ are
the \textit{same} 
if there is an orientation preserving homeomorphism $f$ of $S^3$ 
such that $f(F_1) =F_2$ and $f(K_1) =K_2$; 
otherwise, they are \textit{distinct}. 
It is natural to ask whether
a Seifert fibered surgery $(K, m)$ can have distinct primitive/Seifert positions. 
Recently Guntel \cite{Guntel} has given an infinite family of
such examples.
Her examples are twisted torus knots studied by Dean \cite{Dean}.
Among them, she finds infinitely many pairs of knots
$K_1, K_2$
which have primitive/Seifert positions with the same surface slopes,
and shows that $K_1, K_2$ are actually the same as knots in $S^3$,
but their primitive/Seifert positions are distinct. 

\begin{THM}[\cite{Guntel}]
\label{distinctPS}
There exist infinitely many Seifert fibered surgeries
each of which has distinct primitive/Seifert positions. 
\end{THM}

\begin{REM}
\label{volume}
\textup{In Theorem~$\ref{distinctPS}$, 
we can choose a Seifert fibered surgery $(K, m)$ with distinct primitive/Seifert positions so that 
$K$ is a hyperbolic knot whose complement $S^3 - K$ has an arbitrarily large volume. }
\end{REM}

In the present paper, 
we give yet other families of Seifert fibered surgeries with distinct primitive/Seifert positions 
from a different point of view. 
Our examples are twisted torus knots studied
in \cite{MM3, MM7} (Theorem~\ref{distinctPS1}),
and also Seifert fibered surgeries constructed by the Montesinos trick in \cite{EM1, EM2} (Theorem~\ref{distinctPS2}).
We find infinitely many knots such that
each knot $K$ lies in two genus 2 Heegaard surfaces $F_1, F_2$
with the same surface slopes $m$,
and $(F_1, K, m)$ and $(F_2, K, m)$ are
distinct primitive/Seifert positions.

We use Lemma~\ref{invariant} to show
that two primitive/Seifert positions are distinct.

\begin{LEM}
\label{invariant}
Two primitive/Seifert positions $(F_1, K_1, m)$ and $(F_2, K_2, m)$
for a Seifert fibered surgery $(K, m)$ are distinct
if $i(F_1, K_1, m) \ne i(F_2, K_2, m)$.
\end{LEM}

\textit{Proof of Lemma \ref{invariant}.}
Let us denote the Heegaard splitting of $S^3$ given by 
$F_1$ (resp.\ $F_2$) by $V \cup _{F_1} W$ (resp.\ $V' \cup _{F_2} W'$). 
We may assume that $V[K_1]$ (resp.\ $V'[K_2]$) is a solid torus, 
and $W[K_1]$ (resp.\ $W'[K_2]$) is a Seifert fiber space
with the base orbifold $D^2(p, q)$ (resp.\ $D^2(p', q')$). 
Suppose for a contradiction that we have an orientation preserving homeomorphism $f$ of $S^3$ 
such that $f(K_1) = K_2$ and $f(F_1) = F_2$. 
Then there are two cases to consider: 
$f(V) = V'$ or $f(V) = W'$. 
In the former case $f(W) = W'$ and we have also an orientation preserving homeomorphism
$f_W : W[K_1] \to W'[K_2]$. 
This then implies that $\{p, q\} = \{ p', q' \}$, 
i.e.\ $i(F_1, K_1, m) = i(F_2, K_2, m)$. 
This is a contradiction.
In the latter $f(W) = V'$ and
we have an orientation preserving homeomorphism $f_V : V[K_1] \to W'[K_2]$. 
However, this is impossible
because $V[K_1]$ is a solid torus and $W'[K_2]$
is a Seifert fiber space over the base orbifold
$D^2(p', q')$ $(p', q' \ge 2)$.
\hspace*{\fill} \qed (Lemma~\ref{invariant})

\section{Seifert fibered surgeries which have distinct primitive/Seifert positions I}
\label{distinct1}

\begin{figure}[htbp]
\begin{center}
\includegraphics[width=0.3\linewidth]{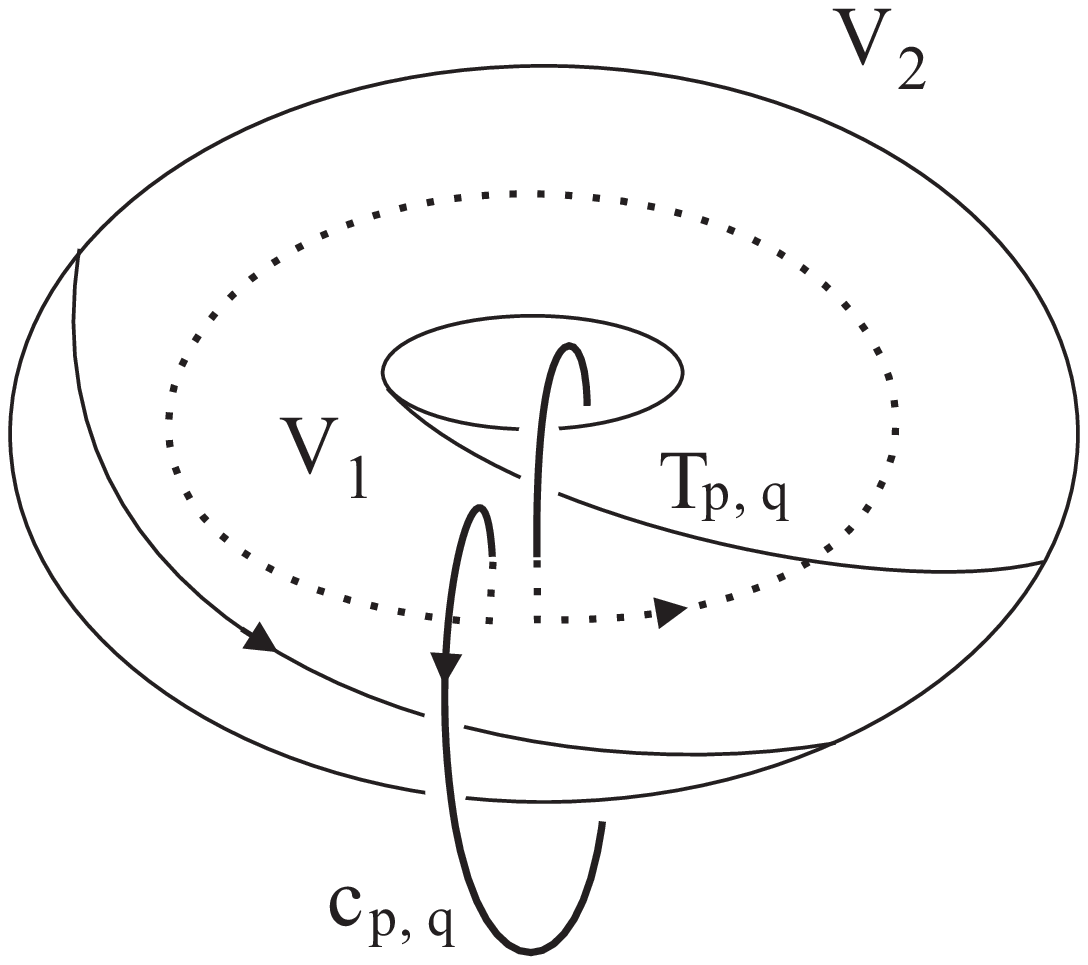}
\caption{}
\label{Tpqcpq}
\end{center}
\end{figure}

Let $V_1$ be a standardly embedded solid torus in $S^3$; 
denote the solid torus $S^3 - \mathrm{int}V_1$ by $V_2$. 
Let $T_{p, q}$ be a torus knot which lies in $\partial V_1$ and 
wraps $p$ times meridionally and $q$ times longitudinally 
in $V_1$. 
Take a trivial knot $c_{p, q}$ in $S^3 - T_{p, q}$
as in Figure~\ref{Tpqcpq}; 
$c_{p, q} \cap V_i$ consists of a single properly embedded arc in $V_i$ 
which is parallel to $\partial V_i$. 
Note that the linking number $\mathrm{lk}(T_{p, q}, c_{p, q})$ 
with orientations indicated in Figure~\ref{Tpqcpq} is $p+q$,
and that $c_{p, q}$ is a meridian of $T_{p, q}$ if $|p+q| = 1$. 
So in the following we assume $|p+q| > 1$. 
We denote by $K(p, q, p+q, n)$ the twisted torus knot
obtained from $T_{p, q}$ by twisting $n$ times along $c_{p, q}$. 
As shown in  \cite[Claim 9.2]{MM1} (\cite[Theorem 3.19(3)]{DMM}), 
$T_{p, q} \cup c_{p, q}$ is a hyperbolic link in $S^3$. 
Hence by \cite[Proposition 5.11]{DMM} 
$K(p, q, p+q, n)$ is a hyperbolic knot if $|n| > 3$. 
In the following, 
for simplicity, 
we denote $c_{p, q}$ by $c$. 

In \cite{MM7} it is shown that 
$(pq + n(p+q)^2)$--surgery on $K(p, q, p+q, n)$
yields a Seifert fiber space 
over $S^2$ with at most three exceptional fibers of indices 
$|p|, |q|, |n|$. 
If $n = 0$, then it is a connected sum of two lens spaces,  
if $n = \pm 1$, then it is a lens space. 
In fact, as shown in \cite{DMM}, 
$(K(p, q, p+q, \varepsilon), pq + \varepsilon (p+q)^2)$
is a Berge's lens surgery \cite{Berge} of Type VII or VIII according as 
$\varepsilon = 1$ or $-1$. 

\begin{THM}
\label{distinctPS1}
Each Seifert fibered surgery $(K(p, q, p+q, n),\ pq + n(p+q)^2)$ 
$(n \ne 0, \pm 1)$ has distinct primitive/Seifert positions. 
\end{THM}

The proof of Corollary~4.8 in \cite{DMM} shows that
for any $r$ there are $p$ and $q$
such that  for infinitely many $n$,
$K(p, q, p+q, n)$ is a hyperbolic knot whose complement
in $S^3$ has volume greater than $r$.
Hence, Theorem~\ref{distinctPS1} implies Theorem~\ref{distinctPS}
and Remark~\ref{volume}.

\begin{figure}[htbp]
\begin{center}
\includegraphics[width=0.5\linewidth]{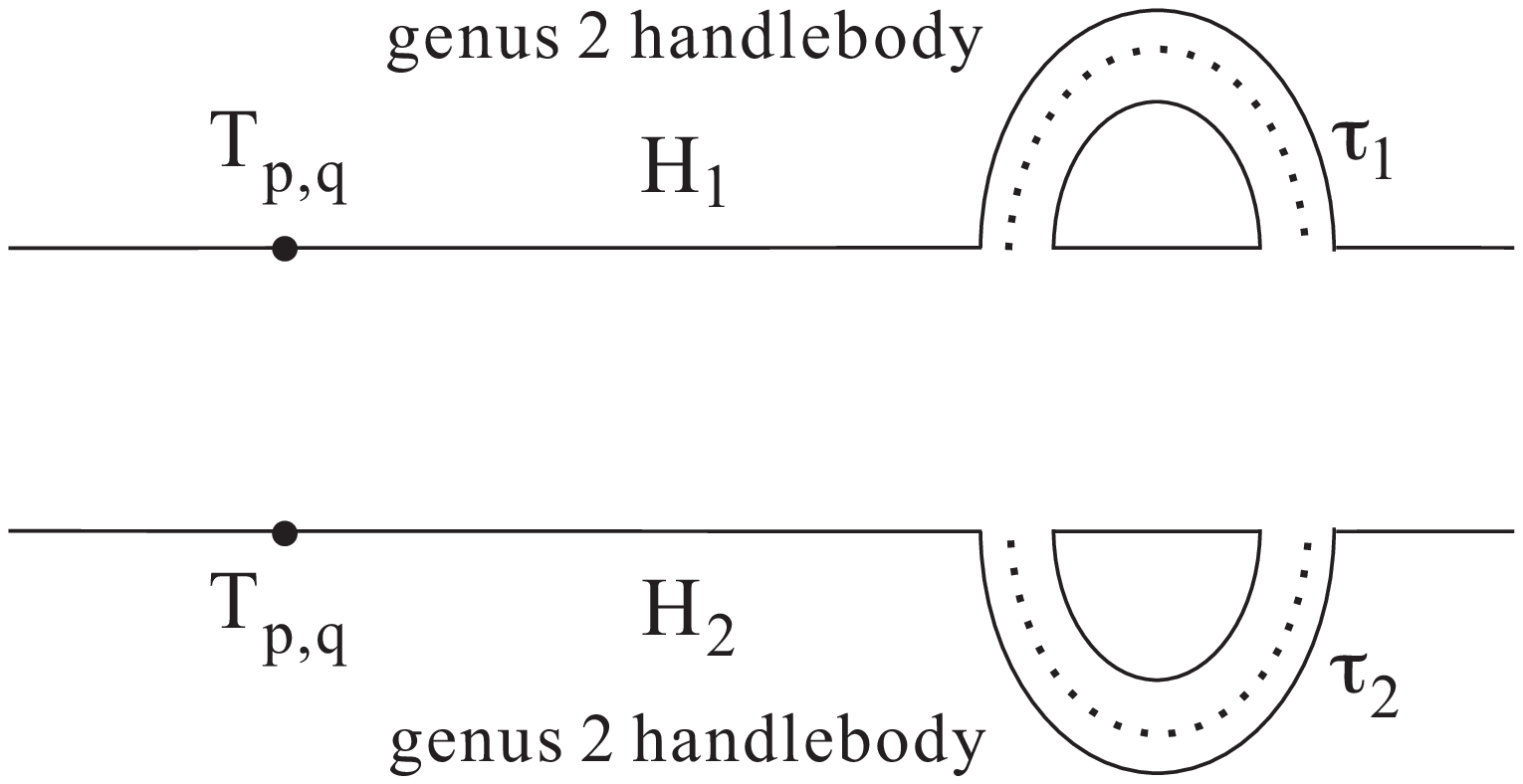}
\caption{}
\label{primitiveSeifert1}
\end{center}
\end{figure}

\textit{Proof of Theorem~\ref{distinctPS1}.}
We follow the argument given in the proof of
\cite[Proposition~5.2]{MM7}. 
Let us put $\tau_i = c \cap V_i$ $(i = 1, 2)$; 
$c = \tau_1 \cup \tau_2$. 
Then 
$H_1 = V_1 - {\rm int}N(c)$ and 
$H_2 = V_2 - {\rm int}N(c)$ are genus $2$ handlebodies, 
and $T_{p, q}$ lies on $\partial V_i - {\rm int}N(c) = H_1\cap H_2$. 
Note that $H_1 \cup H_2 = S^3 - \mathrm{int}N(c)$;
see Figure~\ref{primitiveSeifert1}.

We denote by $U$ the solid torus glued to $S^3 -\mathrm{int}N(c)$
to construct the surgered manifold $c(-\frac{1}{n})$.
Let $A_i = U\cap H_i$; $A_i (\subset \partial H_i)$ is an annulus
whose core is a meridian of $c$.
Since $\tau_i (\subset V_i)$ is parallel to $\partial V_i$,
there is a disk $\Delta_i$ in $H_i$ such that
$\partial \Delta_i$ is the union of an arc in the annulus $A_i$ and
an arc in $\partial H_i -\mathrm{int}A_i$.
Note that $N(\Delta_i) \cup U$ and the closure of $H_i -N(\Delta_i)$
are solid tori,
and their intersection is a disk.
This implies that $H_i \cup U = (N(\Delta_i) \cup U) \cup (H_i -N(\Delta_i))$
is a genus $2$ handlebody for $i =1, 2$.

\begin{LEM}
\label{HG splittings}
$(H_1 \cup U) \cup_F H_2$ and $H_1 \cup_{F'} (H_2 \cup U)$ 
are both genus $2$ Heegaard splitting of $S^3 = c(-\frac{1}{n})$, 
where $F = \partial (H_1 \cup U) = \partial H_2$ 
and $F' = \partial (H_2 \cup U) = \partial H_1$. 
\end{LEM}

Let $\{\mu, \lambda\}$ be a meridian-longitude basis
for $H_1(\partial N(c))$.
Then, a meridian and thus a longitude of $U$ represent
$-n \lambda +\mu$ and $\lambda$ in $H_1(\partial N(c))$,
respectively.
It follows that a meridian of $N(c)$ winds $U$ $n$ times longitudinally.
We thus have the following.

\begin{LEM}
\label{core of A_i}
The core of the annulus $A_i( \subset \partial U)$ winds $U$
$n$ times longitudinally.
\end{LEM}

\begin{figure}[htbp]
\begin{center}
\includegraphics[width=1.0\linewidth]{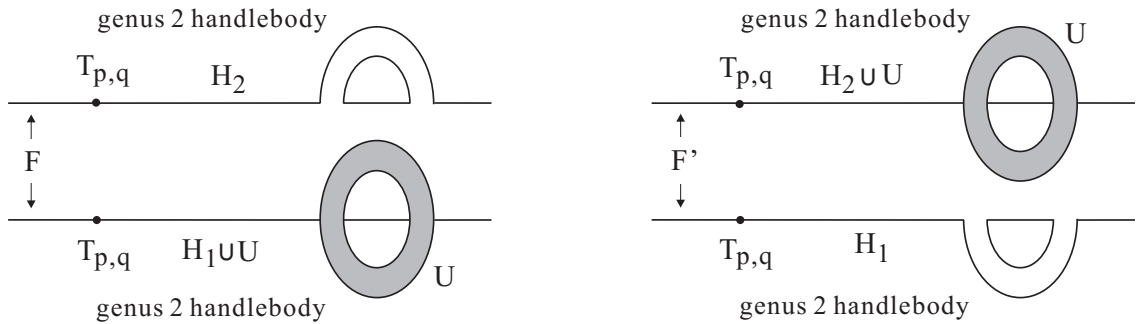}
\caption{Genus $2$ Heegaard splittings carrying $T_{p, q}$}
\label{primitiveSeifert2}
\end{center}
\end{figure}

The twisted torus knot $K(p, q, p+q, n)$ lies on 
$F$ and $F'$. 
See Figure~\ref{primitiveSeifert2}. 
In either case,
the surface slope of $T_{p,q} = K(p, q, p+q, 0)$ is $pq$ and
the surface slope of $K(p, q, p+q, n)$ is the image of that of 
$T_{p,q}$ under $n$-twisting along $c$.
Since $\mathrm{lk}(T_{p,q}, c) = p+q$, 
the surface slope of $K(p, q, p+q, n)$ is $pq + n(p+q)^2$.
\par

\begin{LEM}
\label{arcs in solid tori}
\begin{enumerate}
\item
$H_1[T_{p, q}]$ is a fibered solid torus 
in which the core is an exceptional fiber of index $|q|$ and 
the core of $A_1$ is a regular fiber. 
\item
$H_2 [T_{p, q}]$ is a fibered solid torus 
in which the core is an exceptional fiber of index $|p|$ and 
the core of $A_2$ is a regular fiber.
\end{enumerate}
\end{LEM}

\textit{Proof of Lemma~\ref{arcs in solid tori}.}
See Lemma~9.1 in \cite{MM1}. 
\hspace*{\fill} \qed (Lemma~\ref{arcs in solid tori})

\begin{LEM}
\label{Seifert}
\begin{enumerate}
\item
$(H_1 \cup U)[ K(p, q, p+q, n)]$ is
a Seifert fiber space over $D^2$ with two exceptional 
fibers of indices $|q|$, $|n|$. 
\item
$(H_2 \cup U)[K(p, q, p+q, n)]$ is
a Seifert fiber space over $D^2$ with two exceptional 
fibers of indices $|p|$, $|n|$. 
\end{enumerate}
\end{LEM}

\textit{Proof of Lemma~\ref{Seifert}.}
First observe that $(H_1 \cup U)[K(p, q, p+q, n)]
= H_1[T_{p, q}] \cup U$. 
Since a regular fiber of $H_1[T_{p,q}]$ contained in $A_1$
winds $U$ $n$ times longitudinally
by Lemmas~\ref{core of A_i} and \ref{arcs in solid tori}(1),
$H_1[T_{p, q}] \cup U$ is a Seifert fiber space over
$D^2$ with two exceptional fibers of indices $|q|, |n|$
as claimed in assertion~$(1)$. 
Assertion $(2)$ follows in a similar fashion. 
\hspace*{\fill} \qed (Lemma~\ref{Seifert})

\begin{figure}[htbp]
\begin{center}
\includegraphics[width=1.0\linewidth]{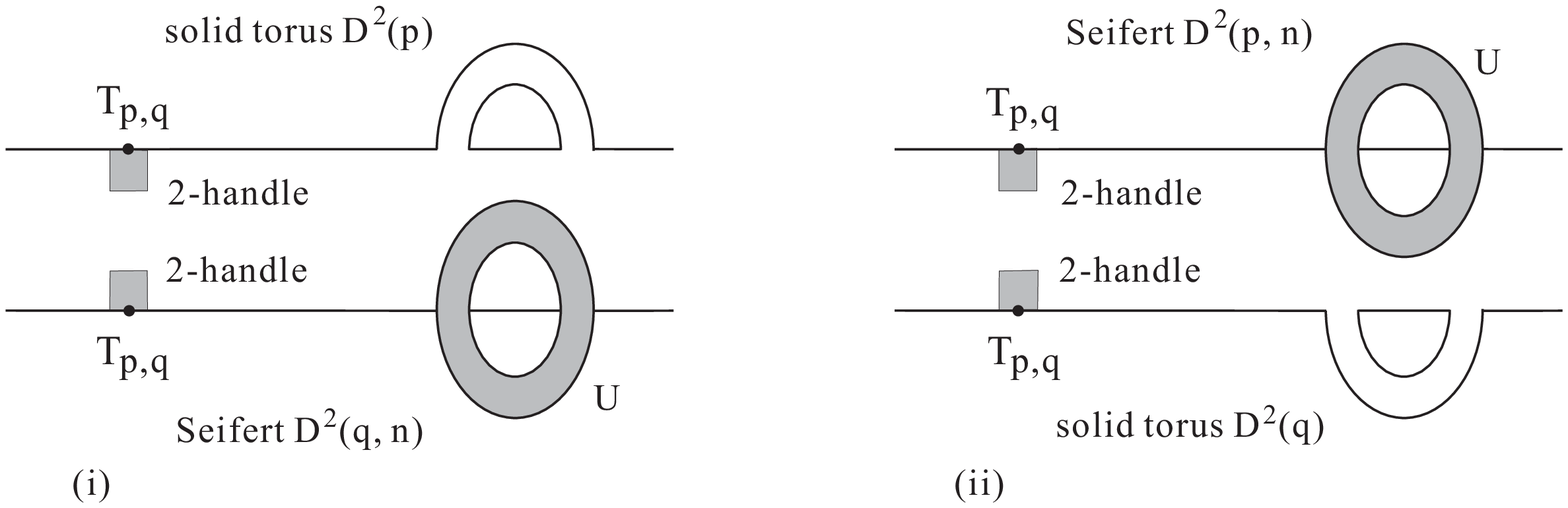}
\caption{}
\label{primitiveSeifert4}
\end{center}
\end{figure}

Therefore the Seifert fibered surgery
$(K(p, q, p+q, n), pq + n(p+q)^2)$ has primitive/Seifert positions
in two ways.
\begin{enumerate}
\item 
$K(p, q, p+q, n)$ is primitive with respect to $H_2$
and Seifert with respect to $H_1 \cup U$;
see Figure~\ref{primitiveSeifert4}(i).
\item
$K(p, q, p+q, n)$ is primitive with respect to $H_1$
and Seifert with respect to $H_2 \cup U$;
see Figure~\ref{primitiveSeifert4}(ii). 
\end{enumerate}

To complete the proof of Theorem~\ref{distinctPS1}
let us show that
$(F, K(p, q, p+q, n), pq + n(p+q)^2)$ and
$(F', K(p, q, p+q, n), pq + n(p+q)^2 )$ are distinct
primitive/Seifert positions.
The former has index $\{|q|, |n|\}$
by Lemma~\ref{Seifert}(1),
and the latter has index $\{|p|, |n|\}$
by Lemma~\ref{Seifert}(2). 
Since $p$ and $q$ are relatively prime, $|p| \ne |q|$.
Then, by Lemma~\ref{invariant} they are distinct.
\hspace*{\fill} \qed (Theorem~\ref{distinctPS1})

\begin{REM}
\label{infinitely many Seifert and Dean's}
\textup{
Twisted torus knots $K(p, q, r, n)$ are obtained
from torus knots $T_{p, q}$, roughly speaking,
by twisting $r$ strands of $T_{p, q}$ $n$ times.
In \cite[Theorem~4.1]{Dean},
Dean obtains five classes of $K(p, q, r, \pm 1)$,
where $p \ge 0, q \ge 0, 0 \le r \le p+q$,
with primitive/Seifert positions.
Let $K$ be a knot in the classes.
We can define $H_1, H_2, U$ for $K$
as for $K(p, q, p+q, n)$ above;
then $K$ is contained in
two genus $2$ Heegaard surfaces $F= \partial(H_1 \cup U),
F' = \partial H_1$.
Although $K(p, q, p+q, n)$ is primitive/Seifert
with respect to both $F$ and $F'$,
in general $K$ is primitive/Seifert with respect to $F$ only. 
}
\end{REM}

\section{Seifert fibered surgeries which have distinct primitive/Seifert positions II}
\label{distinct2}

\begin{figure}[htbp]
\begin{center}
\includegraphics[width=1.0\linewidth]{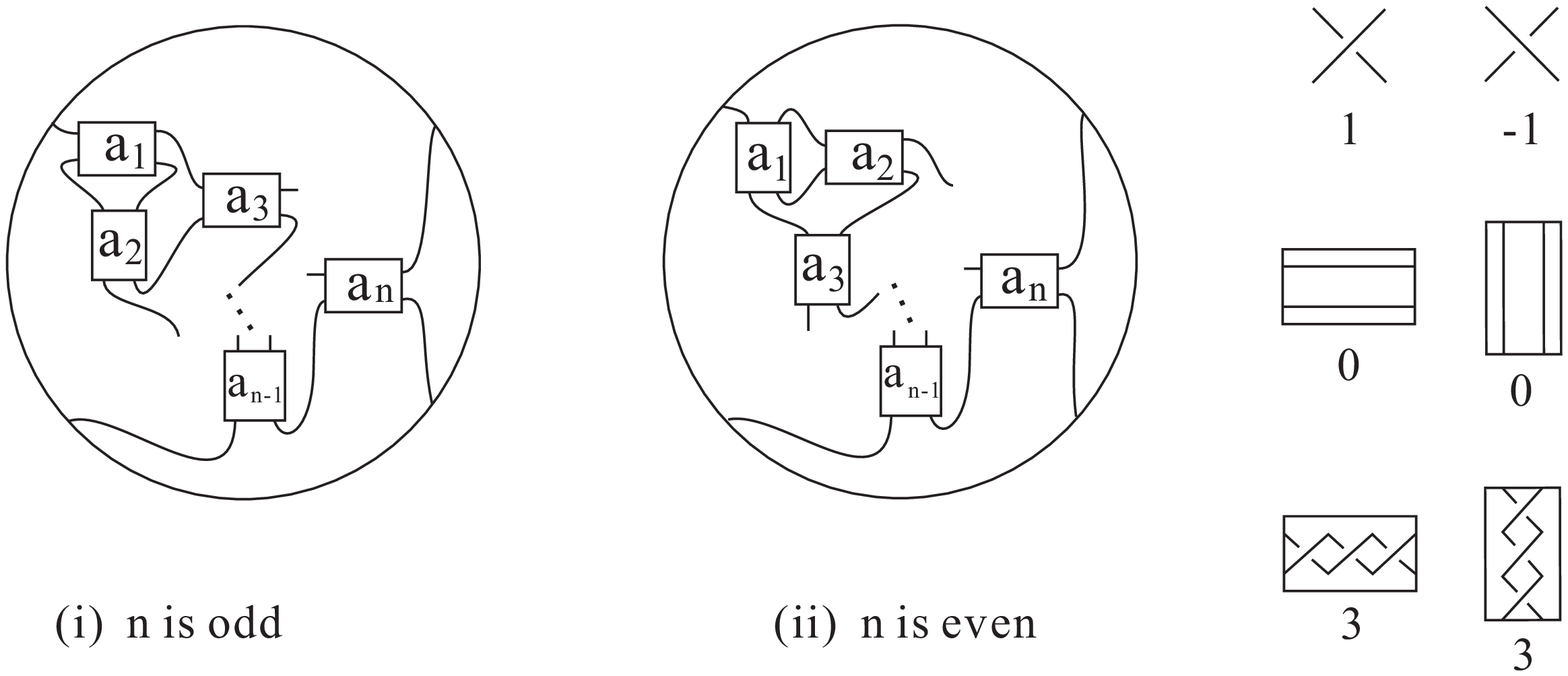}
\caption{Rational tangles}
\label{rtangle}
\end{center}
\end{figure}

A \textit{tangle} $(B, t)$ is a pair of a $3$--ball $B$ and two disjoint arcs $t$ properly embedded in $B$. 
A tangle $(B, t)$ is a \textit{rational tangle} if there is 
a pairwise homeomorphism from $(B, t)$ to the trivial tangle 
$(D^2 \times [0, 1], \{x_1, x_2 \} \times [0, 1])$ 
where $D^2$ is the unit disk and $x_1$ and $x_2$ are 
distinct points in $\mathrm{int}D^2$.
Two rational tangles $(B, t)$ and $(B, t')$ are \textit{equivalent} 
if there is a pairwise homeomorphism  
$h : (B, t) \to (B, t')$ such that $h|_{\partial B} =$ id. 
We can construct rational tangles from sequences of integers 
$[a_1, a_2, \dots, a_n]$ as shown in Figure~\ref{rtangle}.
Denote by $R(a_1, a_2, \dots, a_n)$ the associated rational tangle. 
Each rational tangle can be parametrized by 
$r \in \mathbb{Q} \cup \{ \infty\}$, 
where the rational number $r$ is given by
the continued fraction below. 
Thus we denote the rational tangle corresponding to $r$ by $R(r)$. 

$$ r\ =\ a_{n} + \cfrac 1 {a_{n-1}+ \cfrac 1 {
               \begin{array}{clr}
               \ & & \\[-20pt]
               \ddots & & \\[-17pt]
                      & \ \ + \cfrac{1}{a_1}
               \end{array}
               }}$$
               
\begin{figure}[htbp]
\begin{center}
\includegraphics[width=0.8\linewidth]{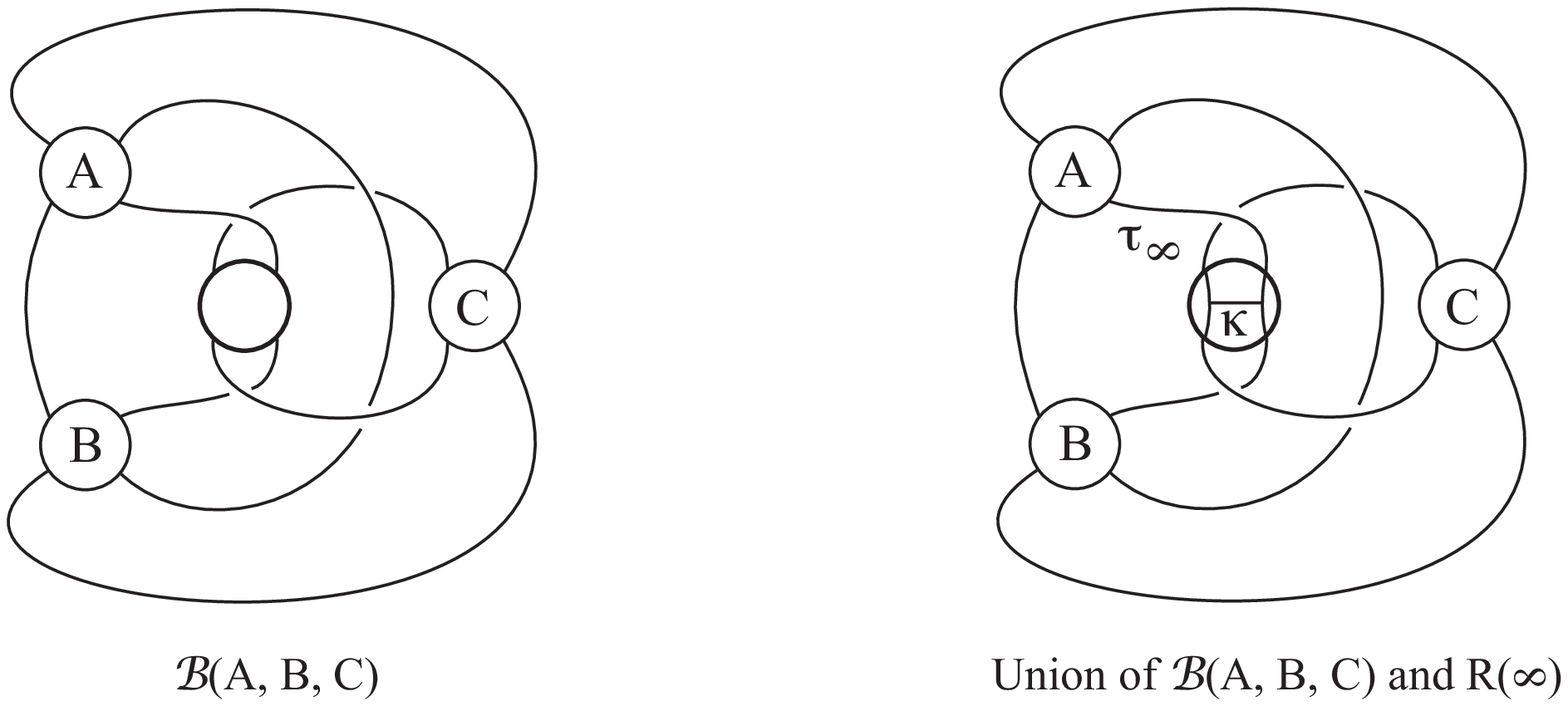}
\caption{}
\label{TABC}
\end{center}
\end{figure}

Let us consider the tangle $\mathcal{B}(A, B, C)$
given by Figure~\ref{TABC}.
In \cite{EM2}, the same tangle $\mathcal{B}(A, B, C)$ is defined
by Figure~9(a) in \cite{EM2}.
However, the figure contains errors;
four crossings of Figure~9(a) in \cite{EM2} should be reversed.
Figure~\ref{TABC} is the corrected diagram.
The union of the  tangle $\mathcal{B}(A, B, C) = (B_1, t_1)$ and a rational tangle $R(s) = (B_2, t_2)$ 
gives a pair $(S^3, \tau_{s}) = (B_1 \cup B_2, t_1 \cup t_2)$. 
We obtain $\tau_s$, a knot or a link in $S^3$. 
In Figure~\ref{TABC} we illustrate the union of $\mathcal{B}(A, B, C)$ and $R(\infty)$. 

In the following,
we assume that $\tau_{\infty}$ is a trivial knot in $S^3$.
Let $\pi_s: \widetilde{S^3}(s) \to S^3$ be 
the two-fold branched covering of $S^3$ along $\tau_s$.
Since $\tau_{\infty}$ is trivial, $\widetilde{S^3}(\infty) =S^3$.
For a subset $X$ of $S^3$,
we often denote $\pi_s^{-1}(X)$ by $\widetilde{X}(s)$,
and $\widetilde{X}(\infty)$ by $\widetilde{X}$ for simplicity.
Let $\kappa$ be an arc
connecting the two vertical strings of $R(\infty)$
as the horizontal arc in Figure~\ref{TABC}.
Then the preimage $\pi_{\infty}^{-1}( \kappa )$ is a knot in $S^3$;
we denote $\pi_{\infty}^{-1}( \kappa )$ by $k(A, B, C)$.
Since the two-fold branched covering of $B_2$
along the rational tangle $t_2$ is a solid torus,
$\widetilde{B_2}(s)$ is a solid torus and in particular
$\widetilde{B_2}$ is a tubular neighborhood of $k(A, B, C)$
in $\widetilde{S^3} =S^3$.
Hence $\widetilde{S^3}(s)$ is obtained from $S^3$ by a Dehn surgery
on $k(A, B, C)$.
We denote the surgery slope by $\gamma_s$.
For $(B_2, t_2) =R(s)$,
if a properly embedded disk $D$ in $B_2 -t_2$ separates the components
of $t_2$, then $\widetilde{D}(s) = \pi_{s}^{-1}( D )$
consists of two meridian disks 
of the glued solid torus $\widetilde{B_2}(s)$.
Hence, a component of $\widetilde{\partial D}(s) = \widetilde{\partial D}$ in $\partial \widetilde{B_2}$
represents the surgery slope $\gamma_s$. \par

Although Figure 9(a) in \cite{EM2} contains errors as mentioned above, 
Lemma~5.1 in \cite{EM2} is correct and we have:

\begin{LEM}[Lemma~5.1 in \cite{EM2}]
\label{trivializable}
$\tau_{\infty}$ is a trivial knot in $S^3$ if either $(1)$ or $(2)$ below holds, 
where $l, m, n, p$ are integers. 
The solutions are the only ones, up to interchanging $A$ and $B$; 
note that there is a rotation interchanging them. 
\begin{enumerate}
\item 
$A = R(l), B = R(m, -l), C = R(-n, 2, m-1, 2, 0)$
\item
$A = R(l), B = R(p, -2, m, -l), C = R(m-1, 2, 0)$
\end{enumerate}
\end{LEM}

In case~$(1)$,
we denote $k(A, B, C)$ by $k(l, m, n, 0)$.  
In case~$(2)$,
we denote $k(A, B, C)$ by 
$k(l, m, 0, p)$. 
As shown in \cite{EM1, EM2}, 
$k(A, B, C)$ are mostly hyperbolic knots. 
See \cite{EM1, EM2} for details.

\begin{figure}[htbp]
\begin{center}
\includegraphics[width=0.8\linewidth]{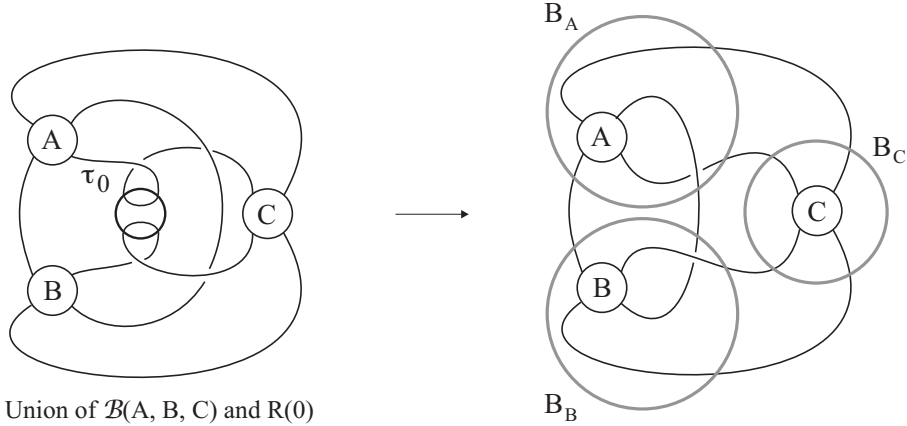}
\caption{$\mathcal{B}(A, B, C) \cup R(0)$}
\label{TABC+0}
\end{center}
\end{figure}

\begin{figure}[htbp]
\begin{center}
\includegraphics[width=1.0\linewidth]{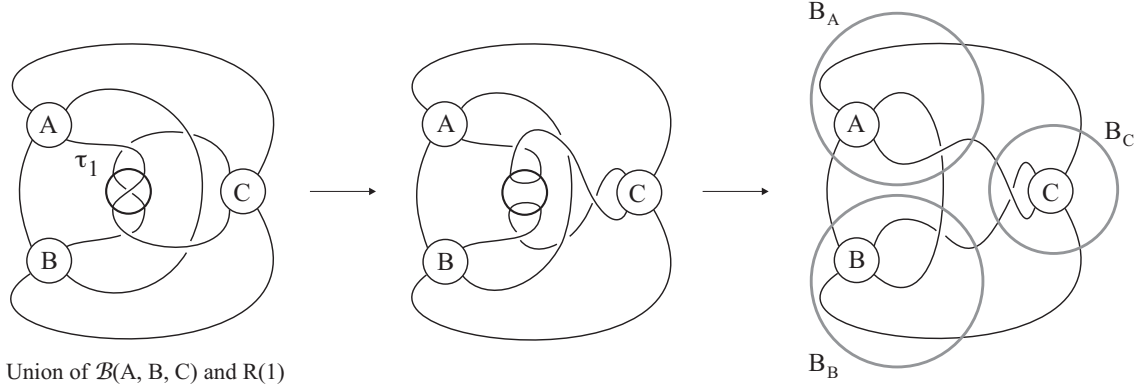}
\caption{$\mathcal{B}(A, B, C) \cup R(1)$}
\label{TABC+1}
\end{center}
\end{figure}

The links $\tau_0$ and $\tau_1$ are Montesinos links with
three branches indicated by the $3$--balls $B_A, B_B, B_C$
in Figures~\ref{TABC+0} and \ref{TABC+1}.
Hence, $\widetilde{S^3}(s) = k(A, B, C)(\gamma_s)$, where $s =0, 1$,
is a Seifert fiber space whose exceptional fibers are the cores 
of $\widetilde{B_A}(s), \widetilde{B_B}(s), \widetilde{B_C}(s)$.
Compute the rational numbers corresponding to
the rational tangles $(B_A, B_A \cap \tau_s),
(B_B, B_B \cap \tau_s), (B_C, B_C \cap \tau_s)$
such that $A, B$, and $C$ satisfy $(1)$ or $(2)$
in Lemma~\ref{trivializable};
then, we obtain the indices of exceptional fibers of
$k(A, B, C)(\gamma_s)$ as follows.
If $(B_X, B_X \cap \tau_s)$ where $X \in \{ A, B, C \}$ corresponds to
a rational number $\frac{p}{q}$, then the Seifert invariant of the core of
$\widetilde{B_X}(s)$ is $-\frac{q}{p}$, and the index is $|p|$. 

\begin{LEM}[corrected Proposition~5.4(2), (3), (5), (6)
in  \cite{EM2}] \
\label{Montesinos}
\begin{enumerate}
\item 
\begin{enumerate}
\item
$\gamma_0$--surgery on $k(l, m, n, 0)$ produces 
a Seifert fiber space over $S^2$ with three exceptional fibers, 
the cores of $\widetilde{B_A}(0), \widetilde{B_B}(0),
\widetilde{B_C}(0)$,
of indices $|l-1|, |l m+m-1|, |2mn -m-n+1|$.
\item
$\gamma_1$--surgery on $k(l, m, n, 0)$ produces 
a Seifert fiber space over $S^2$ with three exceptional fibers, 
the cores of $\widetilde{B_A}(1), \widetilde{B_B}(1),
\widetilde{B_C}(1)$,
of indices $|l+1|,  | lm -m -1 |, |2mn -m +n|$. 
\end{enumerate}

\item
\begin{enumerate}
\item
$\gamma_0$--surgery on $k(l, m, 0, p)$ produces 
a Seifert fiber space over $S^2$ with three exceptional fibers,
the cores of $\widetilde{B_A}(0), \widetilde{B_B}(0),
\widetilde{B_C}(0)$,
of indices $|l-1|, |2l mp - lm- lp + 2mp-m-3p+1|, |m-1|$.
\item
$\gamma_1$--surgery on $k(l, m, 0, p)$ produces 
a Seifert fiber space over $S^2$ with three exceptional fibers, 
the cores of $\widetilde{B_A}(1), \widetilde{B_B}(1),
\widetilde{B_C}(1)$,
of indices $|l+1|, |2lmp -lm -lp -2mp +m -p+1|, |m|$. 
\end{enumerate}
\end{enumerate}
\end{LEM}

\smallskip

In \cite{EM2} a method is given to find primitive/Seifert positions for Seifert fibered surgeries constructed via tangles and double branched covers, 
and this is used to show that each of the surgeries 
$(k(A,B,C),\gamma_s)$ $(s=0,1)$ has a primitive/Seifert position. 
Using this method,  
we prove that each Seifert fibered surgery $(k(A, B, C),\ \gamma_s)$ $(s = 0, 1)$ has distinct primitive/Seifert positions 
if the indices of
the exceptional fibers which are the cores of $\widetilde{B_A}(s)$
and $\widetilde{B_B}(s)$
are not equal. 

\begin{THM}
\label{distinctPS2}
According as $A, B, C$ satisfy
$(1)$ or $(2)$ of Lemma~$\ref{trivializable}$
we assume the following.
\\
If $A, B, C$ satisfy Lemma~$\ref{trivializable}(1)$,
assume that $|l-1| \ne |lm+m-1|$ if $s=0$,
and that $|l+1| \ne | lm -m -1 |$ if $s=1$.
\\
If $A, B, C$ satisfy Lemma~$\ref{trivializable}(2)$,
assume that $|l-1| \ne |2lmp-lm-lp+2mp-m-3p+1|$ if $s= 0$,  
and that $|l+1| \ne |2lmp -lm -lp -2mp +m -p+1|$ if $s=1$. 
\\
Then, each Seifert fibered surgery $(k(A, B, C), \gamma_s)$
$(s = 0, 1)$ has distinct primitive/Seifert positions.
\end{THM}

It follows from \cite[Proposition 5.6]{EM2} that 
the braid index of $k(l, m, n, 0)$ is $2lm - 1$ (resp.\ $2|lm| +1$) 
if $l > 0, m > 0$ (resp.\ $l > 0, m < 0$), 
and that of $k(l, m, 0, p)$ is $2lm -l -1$ (resp.\ $2|lm|+l+1$) 
if $l > 0, m > 0$ (resp.\ $l > 0, m < 0$). 
Hence there are infinitely many knots satisfying the conditions in Theorem \ref{distinctPS2}. 

The assumption in Theorem~\ref{distinctPS2} that the indices of
the exceptional fibers in $\widetilde{B_A}(s)$
and $\widetilde{B_B}(s)$ are not equal is not a necessary condition
for $(k(A, B, C), \gamma_s)$ to have distinct primitive/Seifert positions. Refer to Section~\ref{questions}. 

\begin{figure}[htbp]
\begin{center}
\includegraphics[width=0.8\linewidth]{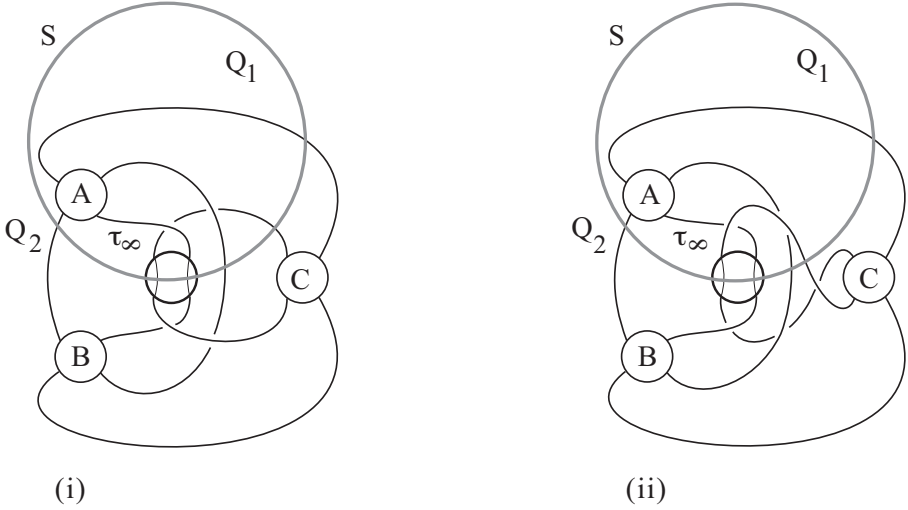}
\caption{$\mathcal{B}(A, B, C) \cup R(\infty)$}
\label{Hsurface1}
\end{center}
\end{figure}

\textit{Proof of Theorem \ref{distinctPS2}.}
Let $s$ be $0$ or $1$.
If $s =0$, let $S$ be the $2$--sphere in $S^3$
shown in Figure~\ref{Hsurface1}(i),
and if $s=1$, let $S$ be the $2$--sphere in $S^3$
shown in Figure~\ref{Hsurface1}(ii).
In either case let $Q_i$ $(i =1, 2)$ be the $3$--balls bounded by $S$
as in Figure~\ref{Hsurface1}.
Note that Figure~\ref{Hsurface1} also describes
the union of the tangles $\mathcal{B}(A, B, C) =(B_1, t_1)$ and
$R(\infty) = (B_2, t_2)$, and $t_1 \cup t_2 = \tau_{\infty}$.
However, $\tau_{\infty}$ in Figure~\ref{Hsurface1}(ii)
is obtained by turning back a portion of $\tau_{\infty}$
in Figure~\ref{TABC}.
The tangles $(Q_i, Q_i \cap \tau_{\infty})$ $(i =1, 2)$
are $3$--string trivial tangles.
Hence, the two-fold branched covering
$\widetilde{Q_1} \cup \widetilde{Q_2}$ gives
a genus $2$ Heegaard splitting of $S^3 =\widetilde{S^3}$,
and $\widetilde{S} = \widetilde{Q_1} \cap \widetilde{Q_2}$ is
a genus $2$ Heegaard surface.
Note that $S \cap B_2$ is a disk intersecting $t_2$ transversely
in two points and containing the arc $\kappa$.
This implies that the annulus $\widetilde{S \cap B_2}$
is a tubular neighborhood of
the knot $k(A, B, C) = \widetilde{\kappa}$
in the Heegaard surface $\widetilde{S}$.
Hence, a component of $\widetilde{S \cap \partial B_2}$
is a simple closed curve in
$\partial \widetilde{B_2} = \partial N(k(A, B, C))$
representing the surface slope of $k(A, B, C)$ in the
Heegaard surface $\widetilde{S}$.

\begin{figure}[htbp]
\begin{center}
\includegraphics[width=0.8\linewidth]{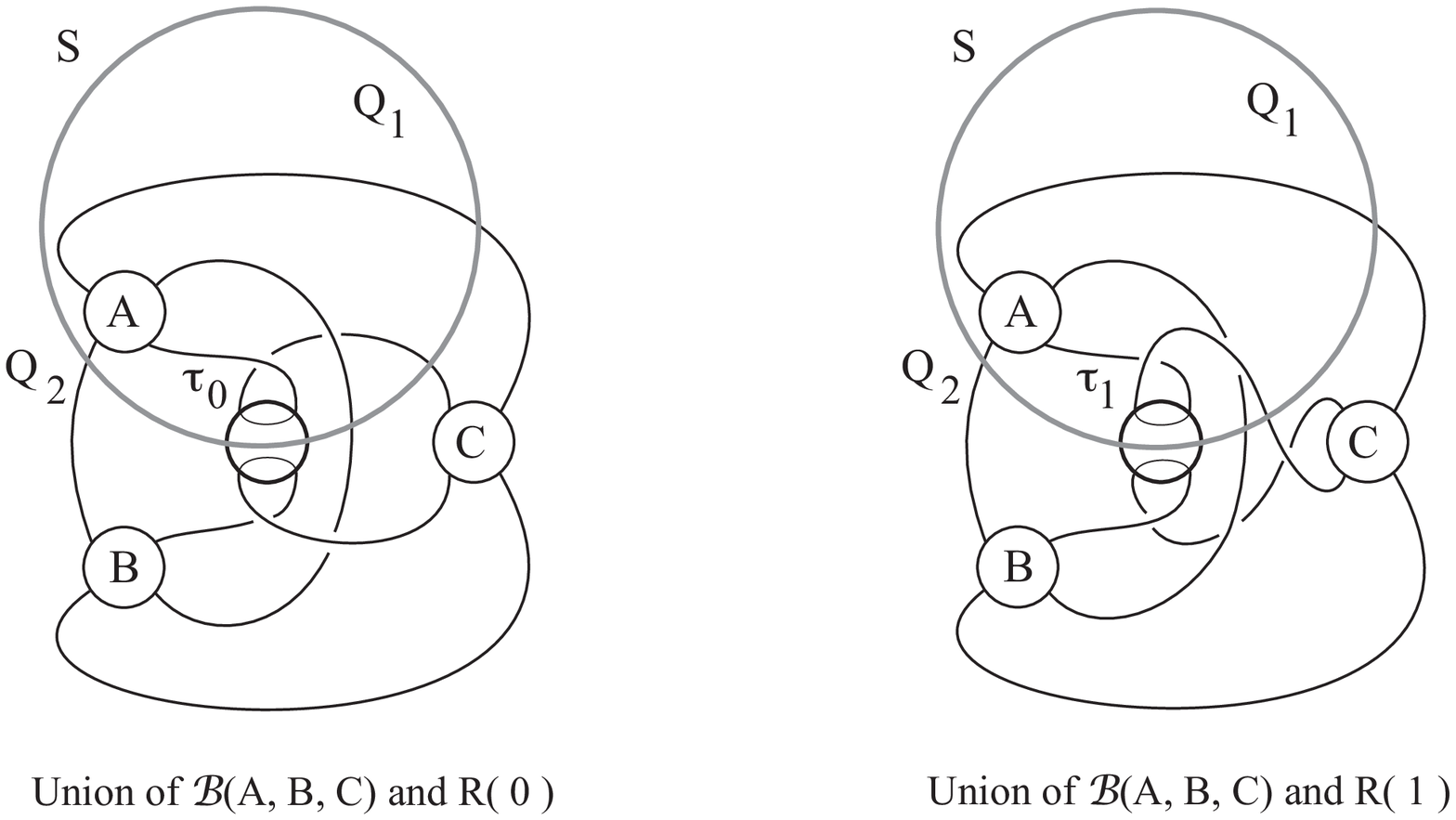}
\caption{}
\label{Hsurface1surgery}
\end{center}
\end{figure}

Now let us show that the surface slope of $k(A, B, C)$
in $\widetilde{S}$ coincides with the surgery slope $\gamma_s$.
Recall that $\gamma_s$--surgery on $k(A, B, C)$
corresponds to replacing $R(\infty)$ with $R(s)$.
The disk $S\cap B_2$ in $S$ is, as shown in Figure~\ref{Hsurface1},
 a ``horizontal'' disk properly embedded in $B_2$.
If $s =0$ and so $S$ is as in Figure~\ref{Hsurface1}(i),
then we may assume that $S\cap B_2$ separates the components of
$t_2$ in $R(0)$ after replaced;
see Figure~\ref{Hsurface1surgery}.
It follows that $\widetilde{S \cap \partial B_2}$ in
$\partial \widetilde{B_2}$ represents the surgery slope $\gamma_s$,
so that the surface slope of $k(A, B, C)$ in $\widetilde{S}$
coincides with $\gamma_0$ as desired.
So assume that $s=1$ and $S$ is as in Figure~\ref{Hsurface1}(ii).
We need to see that the disk $S \cap B_2$ separates the components
of $t_2$ in $R(1)$ attached to $\mathcal{B}(A, B, C)$ in Figures~\ref{TABC}.
The first isotopy in Figure~\ref{TABC+1} turns back
a portion of $\tau_{1}$. 
Then, in the second figure of Figure~\ref{TABC+1},
$t_2$ in $R(1)$ becomes horizontal arcs.
Hence we may assume that the horizontal disk $S\cap B_2$
in Figure~\ref{Hsurface1}(ii) separates the components of $t_2$
in $R(1)$; see Figure~\ref{Hsurface1surgery}.
This implies that a component of $\widetilde{S \cap \partial B_2}$
also represents the surgery slope $\gamma_1$ as desired.

\begin{LEM}
\label{S=ps}
The knot $K = k(A, B, C)$ is in a primitive/Seifert position in
$\widetilde{S}$ with $\gamma_s$ $(s = 0, 1)$ the surface slope,
whose index set is the set of indices of exceptional fibers
in $\widetilde{S^3}(s)$
corresponding to $B_B, B_C$. 
\end{LEM}

\textit{Proof of Lemma \ref{S=ps}.}
We have already shown that the surface slope of $K$ in $\widetilde{S}$
coincides with the surgery slope $\gamma_s$.
We show that $K$ is primitive with respect to
the genus $2$ handlebody $\widetilde{Q_1}$,
and Seifert with respect to $\widetilde{Q_2}$.
First consider $(S^3, \tau_s) =\mathcal{B}(A, B, C) \cup R(s)$.
The $2$--sphere $S$ decomposes $(S^3, \tau_s)$ into
two $2$--string tangles $(Q_1, Q_1 \cap \tau_s)$
and $(Q_2, Q_2 \cap \tau_s)$;
$(Q_1, Q_1 \cap \tau_s)$ is a rational tangle, and
$(Q_2, Q_2 \cap \tau_s)$ is a partial sum of two rational tangles
$(B_B, B_B \cap \tau_s)$ and $(B_C, B_C \cap \tau_s)$
in Figures~\ref{TABC+0}, \ref{TABC+1}.
This implies that $\widetilde{Q_2}(s)$ is a Seifert fiber space
over the disk whose exceptional fibers are the cores of
$\widetilde{B_B}(s), \widetilde{B_C}(s)$.

To complete the proof we prove
$\widetilde{Q_i}(s) \cong \widetilde{Q_i}[K]$.
We consider $(S^3, \tau_{\infty}) =\mathcal{B}(A, B, C) \cup R(\infty)$
again.
The disk $Q_i \cap \partial B_1$ decomposes $Q_i$ into
two $3$--balls $Q_i \cap B_1$ and $Q_i \cap B_2$,
so that $\widetilde{Q_i}(s)
= \widetilde{Q_i \cap B_1}(s) \cup \widetilde{Q_i \cap B_2}(s)$.
Note that $B_1 \cap \tau_s = B_1 \cap \tau_{\infty}$,
and $\tau_s$ intersects $Q_i \cap B_2$ in an arc whose end points
lie in $Q_i \cap \partial B_1$.
Hence, $\widetilde{Q_i \cap B_1}(s) = \widetilde{Q_i \cap B_1}$,
and $\widetilde{Q_i \cap B_2}(s)$ is a $3$--ball attached to
$\widetilde{Q_i \cap B_1}$ along the annulus
$\widetilde{Q_i \cap \partial B_1}$.
In other words,
$\widetilde{Q_i}(s)$ is obtained from $\widetilde{Q_i \cap B_1}$
by attaching a $2$--handle along the annulus
$\widetilde{Q_i \cap \partial B_1}$.
Now replacing $R(s)$ with $R(\infty)$ again,
let us see the relation between $\widetilde{Q_i}$ and
$\widetilde{Q_i \cap B_1}$.
It is not difficult to see the pairwise homeomorphism
$(Q_i \cap B_2, Q_i\cap\partial B_2, Q_i \cap B_2 \cap \tau_{\infty})
\cong
(D^2\times [0, 1], D^2\times\{1\}, \{x_1, x_2\}\times[0, 1])$,
where $x_1, x_2 \in \mathrm{int}D^2$.
This shows that $\widetilde{Q_i \cap \partial B_2}$ is a properly
embedded annulus in $\widetilde{Q_i}$ parallel to
$\widetilde{S \cap B_2}$,
a tubular neighborhood of $K$ in $\widetilde{S}$.
Hence, there is a pairwise homeomorphism from
$(\widetilde{Q_i \cap B_1}, \widetilde{Q_i \cap \partial B_1})$
to $(\widetilde{Q_i}, \widetilde{S \cap B_2})$.
This implies $\widetilde{Q_i}(s) \cong \widetilde{Q_i}[K]$
as desired.
\hspace*{\fill} \qed (Lemma~\ref{S=ps})

\begin{figure}[htbp]
\begin{center}
\includegraphics[width=0.8\linewidth]{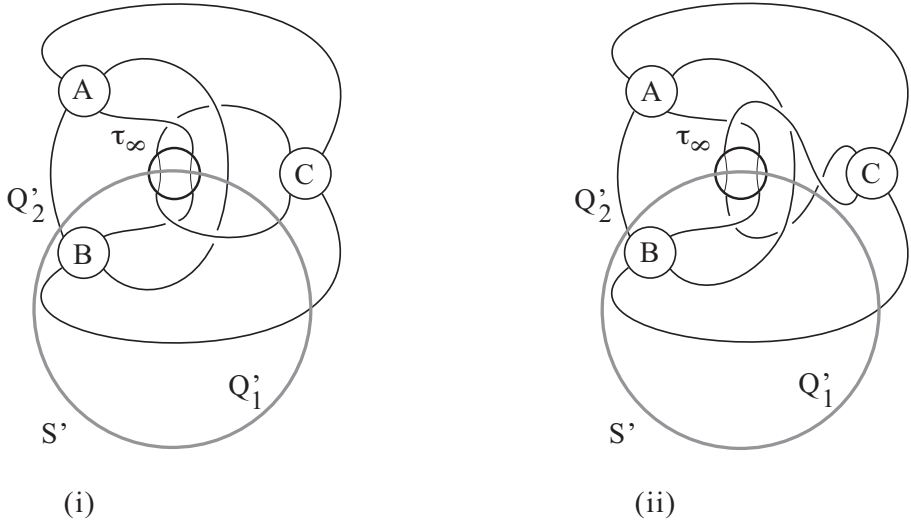}
\caption{$\mathcal{B}(A, B, C) \cup R(\infty)$}
\label{Hsurface2}
\end{center}
\end{figure}

To find yet another primitive/Seifert position of 
$(k(A, B, C), \gamma_s)$, 
take the $2$--sphere $S'$ in $S^3$ as in (i) or (ii) of
Figure~\ref{Hsurface2} according as $s=0$ or $1$.
Let $Q'_i$ $(i =1, 2)$ be the $3$--balls bounded by $S'$
as in Figure~\ref{Hsurface2}.
Then, we can apply the arguments
in the first and second paragraphs
of the proof of Theorem~\ref{distinctPS2}
to $S', Q'_1, Q'_2$ instead of $S, Q_1, Q_2$.
It follows that $\widetilde{Q'_1} \cup
\widetilde{Q'_2}$ is a genus $2$ Heegaard splitting of
$\widetilde{S^3} =S^3$
with $\widetilde{S'}$ the Heegaard surface,
and $k(A, B, C)$ is contained in $\widetilde{S'}$
with $\gamma_s$ $(s =0, 1)$ the surface slope.

We can also apply most of the arguments in the proof of Lemma~\ref{S=ps}.
The only difference is
the fact $(Q'_2, Q'_2\cap \tau_s)$ is a partial sum of
$(B_A, B_A \cap \tau_s)$ and $(B_C, B_C \cap \tau_s)$
instead of $(B_B, B_B \cap \tau_s)$ and $(B_C, B_C \cap \tau_s)$;
see Figure~\ref{Hsurface2surgery}.
Therefore, we obtain Lemma~\ref{S'=ps} below.

\begin{LEM}
\label{S'=ps}
The knot $k(A, B, C)$ is in a primitive/Seifert position in 
$\widetilde{S'}$ with $\gamma_s$ $(s = 0, 1)$ the surface slope,
whose index set is the set of indices of exceptional fibers
in $\widetilde{S^3}(s)$
corresponding to $B_A, B_C$.
\end{LEM}

\begin{figure}[htbp]
\begin{center}
\includegraphics[width=0.8\linewidth]{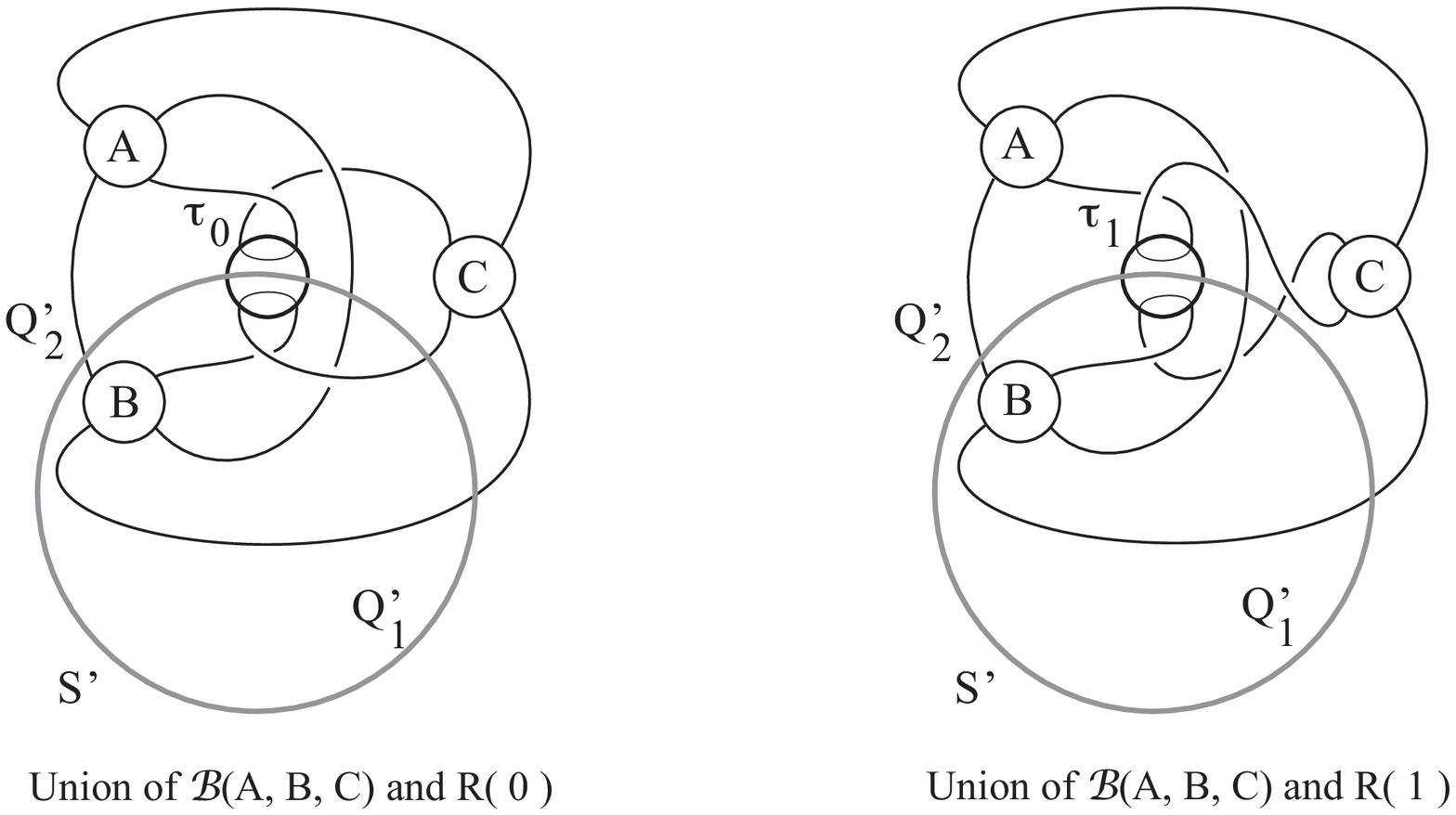}
\caption{}
\label{Hsurface2surgery}
\end{center}
\end{figure}

Recall that by the assumption of Theorem~\ref{distinctPS2}
together with Lemma~\ref{Montesinos},
the indices of the exceptional fibers of
$\widetilde{S^3}(s) = k(A, B, C)(\gamma_s)$
corresponding to $B_A$ and $B_B$
are not equal. 
Lemmas~\ref{S=ps} and \ref{S'=ps} then imply that
the index sets of the primitive/Seifert positions
$(\widetilde{S}, k(A, B, C), \gamma_s)$ and
$(\widetilde{S'}, k(A, B, C), \gamma_s)$ are not equal.
Hence, by Lemma~\ref{invariant}
these are distinct primitive/Seifert positions for
$(k(A, B, C), \gamma_s)$. 
This completes the proof of Theorem~\ref{distinctPS2}. 
\hspace*{\fill} \qed (Theorem~\ref{distinctPS2})

\begin{REM}
\label{infinitely many Seifert}
\textup{
It follows from Lemma \ref{Montesinos} that 
the set $\{ k(A, B, C)(\gamma_s) \}$ $(s = 0, 1)$ consists of infinitely many Seifert fiber spaces. 
If two Seifert fibered surgeries $(k(A, B, C), \gamma_s)$ and $(k(A', B', C'), \gamma'_s)$ are the same, 
then $k(A, B, C)(\gamma_s)$  and $k(A', B', C')(\gamma'_s)$ are homeomorphic. 
Thus the set $\{ (k(A, B, C), \gamma_s) \}$ contains infinitely many Seifert fibered surgeries.
}\end{REM}

\section{Questions}
\label{questions}

In Theorems~\ref{distinctPS1} and \ref{distinctPS2}, 
Seifert fibered surgeries with
distinct primitive/Seifert positions 
have distinct index sets.
However, this is not always the case.
Consider the Seifert fibered surgery $(k(2, 4, n, 0), \gamma_1)$
in Section~\ref{distinct2}. 
The result of $\gamma_1$--surgery on $K = k(2, 4, n, 0)$ is a Seifert fiber space with the base orbifold 
$S^2(\frac{-1}{3}, \frac{4}{3}, \frac{16n-7}{9n-4})$. 
Following Lemma~\ref{S=ps}, we see that 
$(K, \gamma_1)$ has a primitive/Seifert position 
$(\widetilde{S}, K, \gamma_1)$ such that
$\widetilde{Q_2}[K]$ is a Seifert fiber space over the disk 
with Seifert invariants $\frac{4}{3}, \frac{16n-7}{9n-4}$,
where $\widetilde{Q_2}$ is a genus $2$ handlebody
bounded by $\widetilde{S}$. 
Similarly, Lemma~\ref{S'=ps} shows that 
$(K, \gamma_1)$ has a primitive/Seifert position 
$(\widetilde{S'}, K, \gamma_1)$ such that
$\widetilde{Q'_2}[K]$ is a Seifert fiber space over the disk 
with Seifert invariants $\frac{-1}{3}, \frac{16n-7}{9n-4}$,
where $\widetilde{Q'_2}$ is a genus $2$ handlebody
bounded by $\widetilde{S'}$.  
Thus $i(\widetilde{S}, K, \gamma_1) = i(\widetilde{S'}, K, \gamma_1) = \{ 3, |9n-4| \}$. 
If $(\widetilde{S}, K, \gamma_1)$ and 
$(\widetilde{S'}, K, \gamma_1)$ were the same,
then following the argument in the proof of Lemma \ref{invariant}, 
we would have an orientation preserving homeomorphism from 
$\widetilde{Q_2}[K]$ to $\widetilde{Q'_2}[K]$; 
by \cite[VI.18.Theorem]{Ja} the homeomorphism is fiber preserving up to isotopy.
However, since $\frac{4}{3} \not \equiv \frac{-1}{3}$ mod $1$,  
there is no such a homeomorphism \cite[Proposition 2.1]{Hat}.  
Hence the primitive/Seifert positions $(\widetilde{S}, k(2, 4, n, 0), \gamma_1)$ and 
$(\widetilde{S'}, k(2, 4, n, 0), \gamma_1)$ are distinct.

\begin{QUES}
\label{q1}
\textup{
Does there exist a Seifert fibered surgery
which has distinct primitive/Seifert positions 
$(F_1, K_1, m)$ and $(F_2, K_2, m)$ satisfying the following condition?
\\
Condition.
Let $W_i$ $(i =1, 2)$ be a genus $2$ handlebody bounded by $F_i$
with respect to which $K_i (\subset \partial W_i)$ is Seifert.
Then there is an orientation preserving homeomorphism
from $W_1[K_1]$ to $W_2[K_2]$.} 
\end{QUES}

Even if a Seifert fibered surgery $(K, m)$ has distinct primitive/Seifert positions, 
we expect that the number of such positions is not so large. 
In fact,
we do not even have an example of a Seifert fibered surgery which has three primitive/Seifert positions. 

\begin{QUES}
\label{q2}
\textup{
Does there exist a universal bound for the number of primitive/Seifert positions for 
a Seifert fibered surgery? }
\end{QUES}

\textbf{Acknowledgments.}\ 
We would like to thank the referee for careful reading and useful suggestions. 
The first author was partially supported by PAPIIT-UNAM grant IN102808.
The third author has been partially supported by JSPS Grants--in--Aid for Scientific 
Research (C) (No.21540098), 
The Ministry of Education, Culture, Sports, Science and Technology, Japan.


\begin{thebibliography}{99}

\bibitem{Berge} J. Berge; 
Some knots with surgeries yielding lens spaces, 
unpublished manuscript.

\bibitem{Dean} J. Dean; 
Small Seifert-fibered {D}ehn surgery on hyperbolic knots, 
Algebr. Geom. Topol.\ \textbf{3} (2003), 435--472

\bibitem{DMM} A. Deruelle, K. Miyazaki and K. Motegi; 
Networking Seifert Surgeries on Knots, 
to appear in Mem.\ Amer.\ Math.\ Soc. 

\bibitem{EM} M. Eudave-Mu\~noz; 
Band sums of links which yield composite links. 
The cabling conjecture for strongly invertible knots, 
Trans.\ Amer.\ Math.\ Soc.\ \textbf{330} (1992), 463--501.  

\bibitem{EM1} M. Eudave-Mu\~noz;
Non-hyperbolic manifolds obtained by Dehn surgery on a hyperbolic knot,
In: Studies in Advanced Mathematics vol.\ \textbf{2}, part~1,
(ed.\ W. Kazez), 1997, Amer.\ Math.\ Soc.\ and International Press,
pp. 35--61.

\bibitem{EM2} M. Eudave-Mu\~noz; 
On hyperbolic knots with Seifert fibered Dehn surgeries, 
Topology Appl.\ \textbf{121} (2002), 
119--141. 

\bibitem{Go90} C.McA. Gordon;  
Dehn surgery on knots, 
Proceedings ICM Kyoto 1990, (1991), 555--590. 

\bibitem{Guntel} B. Guntel;  
Knots with distinct primitive/primitive and primitive/Seifert representatives, 
to appear in J.\ Knot Theory Ramifications. 

\bibitem{Hat} A. E. Hatcher; 
Notes on basic $3$-manifold topology, 
freely available at \texttt{http://www.math.cornell.edu/\~{}hatcher}, 
2000. 

\bibitem{Ja} W. Jaco; 
Lectures on three manifold topology, 
CBMS Regional Conference Series in Math.\ 43, Amer.\ Math.\ 
Soc., 1980. 

\bibitem{MMM} T. Mattman, K. Miyazaki and K. Motegi;
Seifert fibered surgeries which do not arise from primitive/Seifert-fibered constructions, 
Trans. Amer. Math. Soc.\ \textbf{358} (2006), 4045--4055. 

\bibitem{MM1} K. Miyazaki and K. Motegi; 
Seifert fibred manifolds and Dehn surgery, 
Topology \textbf{36} (1997), 579--603. 

\bibitem{MM3} K. Miyazaki and K. Motegi; 
Seifert fibered manifolds and Dehn surgery III, 
Comm.\ Anal.\ Geom.\ \textbf{7} (1999), 551--582. 

\bibitem{MM7} K. Miyazaki and K. Motegi; 
On primitive/Seifert-fibered constructions, 
Math.\ Proc.\ Camb.\ Phil.\ Soc.\ \textbf{138} (2005), 421--435. 

\bibitem{Mon} J. M. Montesinos; 
Surgery on links and double branched coverings of $S^3$, 
Ann.\ Math.\ Studies \textbf{84} (1975), 227--260. 

\bibitem{Tera} M. Teragaito; 
A Seifert fibered manifold with infinitely many knot-surgery descriptions, 
Int. Math. Res. Not.\ \textbf{9} (2007), Art. ID rnm 028, 16 pp.


\bibitem{Wa} F. Waldhausen; 
Heegaard-Zerlegungen der $3$--Sph\"are, 
Topology \textbf{7} (1968), 195--203. 

\end{thebibliography}
\end{document}